\documentclass[sn-mathphys]{sn-jnl}

\jyear{2021}%

\theoremstyle{thmstyleone}%
\newtheorem{theorem}{Theorem}
\newtheorem{proposition}[theorem]{Proposition}%

\usepackage{graphicx}
\usepackage[pdftex]{pict2e}  

\theoremstyle{thmstyletwo}%

\theoremstyle{thmstylethree}%

\newcommand{\Delete}[1]{}
\newcommand{\pend}{\hspace*{\fill} $\Box$}

\begin{document}

\title[Sign representation of single-peaked preferences and Bruhat orders]{Sign representation of single-peaked preferences and Bruhat orders}


\author[1]{\fnm{Ping} \sur{Zhan}\footnote{\ zhan@edogawa-u.ac.jp \par Department of Communication and Business, Edogawa University \par 474 Komaki, Nagareyama, Chiba, 270-0198, Japan}}

\abstract{Single-peaked preferences and domains are extensively researched in social science and economics. 
In this study, we examine the interval property as well as combinatorial structure of single-peaked preferences on a fixed Left-Right social axis.
We introduce a sign representation of single-peaked preferences; consequently, some cardinalities of single-peaked domains are easily obtained.  
Basic operations on the sign representation, which completely define the Bruhat poset, are also provided.
The applications to known results and an isomorphic relation with associated rhombus tiling are given.
Finally, we some discussions of related topics. }




\maketitle

\section{Introduction}\label{sec1}
We consider elections or ballots with two finite sets, candidates and votes. 
Voters' single-peaked preferences are linear orderings defined on the set of candidates or alternatives.

The study of single-peaked preferences was introduced by Black  \cite{black1948} such that the simple aggregation majority voting rule does not yield cycles. 
Single-peakedness preserves strategy-proofness in assignments \cite{moulin1980}.
The structure and characterization of single-peaked domains have been researched extensively  
\cite{puppe2018, EFS2020}, see also a review paper \cite{monjardet2009}.
Recently, single-peakedness is studied in various topics, for applications in matching and assignment \cite{bade2019}, counting and distribution \cite{karpov2020}, construction by tiling \cite{zhan2019}, forbidden configuration \cite{BH2011}, and finding optimal committees under proportional approval voting \cite{peters2018}.

In this study, we introduce an interval and sign representation of single-peaked preferences with a fixed one-dimension axis, as well as related basic operations. 
To the best of our knowledge, the sign representation of single-peakedness is new.
Some applications and discussions are also provided.
Specifically, this study contributes the following.
\begin{itemize}
\item[] 
Based on the interval property of single-peaked linear order, we introduce a sign representation that uniquely specifies each single-peaked preference. 
Some cardinality number related to single-peaked domains are easily obtained using the representation. 
The operations based on the sign representation which specify the structure of Bruhat orders are also provided.
The sign representation and operations can be used to characterize the minimal richness and semi-connectedness of single-peaked domains \cite{puppe2018, EFS2020} in a different way. 
The procedural intervals and sign representation are useful in designing algorithms related to the single-peakedness \cite{bade2019}. 
We expect that will initiate new researches.  
\end{itemize}

The rest of this manuscript is organized as follows. 
Section \ref{sec:def} provides definitions. 
In Section \ref{sec:connect}, after examining some combinatorial structures, we present a sign representation of single-peaked preferences, and the associated operations; then show some applications obtained by these results. 
Finally, concluding remarks are provided in Section \ref{sec:conclud}.
\section{Definitions and Preliminaries} \label{sec:def}
\subsection{Single-peaked Preferences and Single-peaked Domains} \label{sec:def_single}
We consider an election with two finite sets, candidates or alternatives $X \equiv \{x_1,\dots,x_n\}$, votes $V\equiv \{ v_1, \dots , v_m \}$. 
For alternatives $x_i,x_j \in X$, the relation $x_i \succ_\ell x_j$ means that voter $v_\ell$ strictly prefers $x_i$ to $x_j$.
A vote's \emph{preference} is a linear order (i.e., complete, transitive, and asymmetric binary relations) on $X$: in other words, preference $\succ _\ell \ \equiv \ x_1 \succ_\ell \cdots \succ_\ell x_n$ with $x_j \in X$ $(1 \le j \le n
)$. 
We simplify $x_i \succ_\ell x_j$ as $x_i \succ x_j$ if we need not to specify a voter. 
We write $x_ix_j$ instead of $x_i \succ x_j$  (if there is no confusion). 

Let $\mathcal{L}(X)$ denote \emph{all} linear orders (or permutations) on $X$. 
A subset $\mathcal{D} \subseteq \mathcal{L}(X) $ is called a \emph{domain of preferences} or simply a \emph{domain}.
A \emph{profile} $\pi_m=(\succ_1,\cdots,\succ_m)$ on domain ${\cal D}$ is an element of the Cartesian product ${\cal D}^m$ of $m$ voters. 
For simplicity, suppose $m$ being an odd number.

The \emph{$($simple$)$ majority relation} associated with a profile $\pi_m$ is defined by $x_i \prec^{maj}_{\pi_m} x_j $ if (and only if) more voters ($>m/2$) prefer $x_j$ than $x_i$. 
\emph{Condorcet domains} (briefed as CDs) are subsets of $ \mathcal{L}(X)$ with the property that, whenever the preferences of all voters belong to this set, the majority relation has no cycles \cite{fishburn1996, DK2013}. 
For example, let $\pi_3=(x_1x_2 x_3,x_2x_3x_1,x_3x_1x_2)$, it is easy to check that 
\begin{equation}
x_1 \prec^{maj}_{\pi_3} x_3 \prec^{maj}_{\pi_3} x_2 \prec^{maj}_{\pi_3} x_1 ,   \nonumber
\end{equation} 
i.e., a \emph{majority cycle}, $x_i \prec^{maj}_{\pi_m} ,\cdots, \prec^{maj}_{\pi_m} x_i$ for some $x_i \in X$, occurring.  
Hence, $\{x_1x_2 x_3,x_2x_3x_1,x_3x_1x_2\}$ is not a CD.

Given a positive integer $j$, we denote the set $\{1, \dots, j\}$ by  $[j]$.
We consider single-peaked preferences on a fixed Left-Right social axis on alternative $X$, which, up to an isomorphism, is $[n]$ with $1 < 2 < \cdots < n$ (refer to, e.g., \cite{slinko2019} for more general case). 
Let $ k \in [n]$ be the peak, or the top of the preference. 
A preference is called \emph{single-peaked} if for all $k_1, k_2 \in [n]$,  $k_2<k_1\le k$ or $k_2 > k_1 \ge k$ implies  $ k_1 \succ k_2$. 
A domain is single-peaked if each preference in the domain is single-peaked. 

Furthermore, the domain of \emph{all} orders on $[n]$ that are single-peaked is denoted by $\mathcal{SP}([n])$  \cite{puppe2018}. 
Here is an example of a single-peaked domain: $\mathcal{SP}([4]) =\{1234, 2134, 2314, 2341, 3214, 3241, 3421, 4321 \} $. 
It is well-known that single-peaked domains are special cases of CDs. 
\subsection{Bruhat Orders}
Let $\Omega = \{(i,j) \mid i,j \in [n], i <j \}$. 
For a linear order $\sigma$, a pair $(i,j) \in \Omega $ is called an \emph{inversion} w.r.t. $\sigma$ if $ j \prec_{\sigma} i$ and $i<j$. 
The set of all inversions for $\sigma$ is denoted by Inv$(\sigma)$. 
For linear orders $\sigma, \tau \in \mathcal{L}=\mathcal{L}([n])$, we write $\sigma \ll \tau$ if ${\rm Inv}(\sigma) \subseteq {\rm Inv} (\tau)$.
The \emph{Bruhat order} is the partial order on $\{ {\rm Inv}(\sigma) \mid \sigma \in \mathcal{L}([n]) \} $ induced by the set inclusion.
We denote \emph{Bruhat poset} $(\mathcal{D}, \ll)$ on domain $\mathcal{D} \subseteq \mathcal{L}([n])$ by $\mathbb{B}_n(\mathcal{D},1)$, or $\mathbb{B}(\mathcal{D},1)$, also refer to  \cite{abello_j1991, g_reiner08} and \cite{ziegler1993} for more general cases. 
The \emph{Bruhat digraph} is defined by drawing an arc from $\sigma $ to $\tau$ if $\tau$ covers $\sigma$, i.e., Inv$(\sigma) \subset {\rm Inv}(\tau)$ and ${\rm Inv}(\tau)$ has exactly one more inversion than Inv$(\sigma)$ \cite{DKK2012}.  
Fig. \ref{fig:bruhat_def} (a) shows an example of a linear order with $n=3$; Fig. \ref{fig:bruhat_def} (b) gives the Bruhat poset with inversion pairs corresponding to the linear orders in Fig. \ref{fig:bruhat_def} (a). 
\setlength{\unitlength}{1mm}
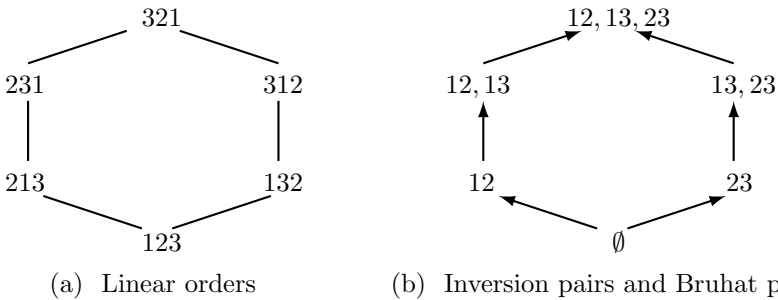
\begin{figure}[h!]
 \centering
\begin{picture}(155,38)
\linethickness{0.6pt}

\put(11,5.5){ \begin{minipage}[t]{50mm}
\begin{picture}(35,31)(0,0)
\linethickness{0.8pt}
\put(2,25){\line(3,1){13}} \put(35,25){\line(-3,1){13}} 
\put(2,12){\line(0,1){8}} \put(35,12){\line(0,1){8}}
\put(21,3){\line(3,1){13}} \put(17,3){\line(-3,1){13}}

\put(17,30){$321$} 
\put(-1,21){$231$}  \put(33,21){$312$} 
\put(-1,8){$213$}  \put(33,8){$132$} 
\put(17,0){$123$}
\end{picture}
\end{minipage} }

\put(71,5.5){ \begin{minipage}[t]{50mm}
\begin{picture}(35,30)(0,0)
\linethickness{0.8pt}
\put(2,25){\vector(3,1){13}} \put(35,25){\vector(-3,1){13}} 
\put(2,12){\vector(0,1){8}} \put(35,12){\vector(0,1){8}}
\put(21,3){\vector(3,1){13}} \put(17,3){\vector(-3,1){13}}

\put(13,30){$12,13,23$} 
\put(-3,21){$12,13$}  \put(32,21){$13,23$} 
\put(0,8){$12$}  \put(34,8){$23$} 
\put(19,0){$\emptyset$}
\end{picture}
\end{minipage} }
\put(17,-0){(a) \ Linear orders}
\put(62,-0){(b) \ Inversion pairs and Bruhat poset}
\end{picture}
\caption{(a) Linear orders, and  (b) the associated Bruhat poset (digraph) with $n=3$} 
\label{fig:bruhat_def}
\end{figure}
\Delete{
A linear order of 
$\Omega $ is \emph{admissible} if the 2-subsets of any 3-subset of $[n]$ appear either in lexicographic or reverse lexicographic (anti-lexicographic) order. 
Precisely, for any $ i < j < k , $ set $\{\{i,j\}, \{i,k\}, \{j,k\}\}$ appears either in $(\{i,j\},\{i,k\}, \{j,k\})$ or in $(\{j,k\}, \{i,k\}, \{i,j\})$ order, see the papers given by  \cite{DKK2021, g_reiner08, ziegler1993} for more detailed and general discussions. 
Note that the above definition slightly differs from the one given by \cite{g_reiner08} since we are interested in the order of the voters on single-peaked and peak-pit domains.
The \emph{higher Bruhat order}, denoted by $\mathbb{B}_n(\mathcal{D},2):= (\mathcal{D},\ll_{adm})$ or $ \mathbb{B}(\mathcal{D},2)$ on domain $\mathcal{D} \subset \mathcal{L}([n])$, is a partial order on admissible sets denoted by $\ll_{adm}$, which is defined by the transitive closure of a single step set inclusion. 
In Fig. \ref{fig:bruhat_def} (b), the left and the right paths are in lexicographic and reverse lexicographic orders, respectively. 
}
\subsection{Rhombus Tilings of Preference Domains} \label{sec:tile_def}
The representation of the preference domains by rhombus tilings was developed by Danilov et al. \cite{DKK2012}. 

Let $\xi_i \in \mathbb{R} \times \mathbb{R}_{>0}$ ($i \in N$) be vectors in the upper half-plane. 
A (graphic) \emph{rhombus tiling} (or simply a \emph{tiling}), denoted by $T$,
is a symmetric \emph{zonogon} formed by $\sum_{i-1}^n \chi_i \xi_i$ over all $\chi_i \in  \{0,1\}$ $1\le i \le n$),  
with the \emph{source} vertex $(0,0)$ at the bottom, and the \emph{sink} vertex $\xi_1+\xi_2+ \cdots +\xi_n$ at the top. 
A rhombus congruent to the sum of two vectors $\xi_i$ and $\xi_j$ is called an $ij$-\emph{tile} ($i \ne j$), or simply a \emph{tile}. 
Borrowing the language of a graph, tiling $T$ is acyclic since all arcs are oriented upward. 
A maximal path from the source to the sink is called a \emph{snake}, which contains exactly one $i \in [n]$. 
The labels on the directed path of each snake form  a linear order, or a permutation of alternatives of $[n]$. 
Also refer to \cite{DKK2021, zhan2019} for more details of rhombus tilings and preference domains, refer to \cite{elnitsky1997} for more general linear orders and permutations.

Fig. \ref{fig:single_tiling_4} is an example of rhombus tiling for $n=4.$
We label the arcs of tiling with associated alternatives $1 \le i \le 4$ (same labels of parallel arcs are omitted if they are obvious, e.g., all green dashed lines are labelled as ``2"). 
In Fig. \ref{fig:single_tiling_4}, all directed paths from the bottom to top form the single-peaked domain $\mathcal{SP}([4])$.

\setlength{\unitlength}{1mm}
\begin{figure}[h!]
 \centering
\begin{picture}(155,44)
\linethickness{0.6pt}

\put(38,4.5){ \begin{minipage}[t]{50mm}
\begin{picture}(35,35)(0,0)
\linethickness{0.5pt}

\put(30.8,31){\vector(-2,1){13}}
\put(18,24.5){\vector(2,1){13}}

\put(0.5,19){\vector(2,5){5}} 
\multiput(35.5,19)(-0.8,2){5}{\textcolor{green}{\line(-2,5){0.45}}}
\put(31.5,29){\vector(-2,5){1}}

\put(13,12.5){\vector(-2,1){13}} 
\put(23,12.5){\vector(2,1){13}} 
\multiput(23,12.5)(-0.8,2){5}{\textcolor{green}{\line(-2,5){0.45}}}
\put(19,22){\vector(-2,5){1}}

\multiput(5.2,6.5)(-0.8,2){5}{\textcolor{green}{\line(-2,5){0.45}}}
\put(1.5,16.5){\vector(-2,5){1}}
\put(30.8,6.5){\vector(2,5){5}}

\put(18,0){\vector(2,1){13}}
\put(18,0){\vector(-2,1){13}}
\put(18,0){\vector(2,5){5}}

\linethickness{1.7pt}
\multiput(18,0)(-0.8,2){5}{\textcolor{green}{\line(-2,5){0.45}}}
\put(14,10){\vector(-2,5){1}}

\put(13,12.5){\vector(2,5){5}} 
\put(18,24.5){\vector(-2,1){13}}
\put(5.2,31){\vector(2,1){13}}

\put(10,35){$4$}  \put(25,35){$1$} 
\put(0,25){$3$}  \put(34,25){$2$} 
\put(0,9){$2$}  \put(34,9){$3$} 
\put(10,0){$1$} \put(25,0){$4$}
\end{picture}
\end{minipage} }
\end{picture}
\caption{ Tiling of a single-peaked domain $\mathcal{SP}([4])$, where the thick (including the dashed) lines form a snake representing the preference of $2314$. } 
\label{fig:single_tiling_4}
\end{figure}
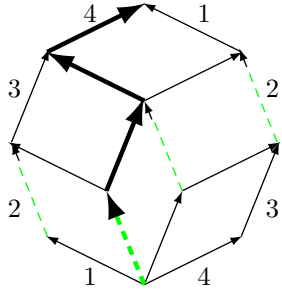

A domain $\mathcal{D} \subset \mathcal{L}([n])$ is called \emph{(single-pit) pit domain} if restricted to any triple $1 \le i < j <k \le n$, we have $\{ijk, ikj, kij, kji \}$.
A domain $\mathcal{D} \subseteq \mathcal{L}([n]) $ is called as a \emph{peak-pit} domain if restricted to $\{i,j,k \} \subseteq [n]$, $\mathcal{D}_{\{i,j,k \}}$ is either a single-peaked, or single-pit domain.  
The fact that each maximal peak-pit domain can be represented by a rhombus tiling has been shown by \cite{DKK2012}. 

In the tiling representation of preference domains, we assume that the left bound is $12\cdots n$ and the right bound is $n(n-1) \cdots 1$. 
There are exactly $\binom{n}{2}= \frac{n(n+1}{2}$ tiles, or reverse pairs, in each tiling on $[n]$. 
The right snakes have more inversion pairs than the left ones in a  Bruhat poset, all snakes of a tiling also form a lattice
 \cite{elnitsky1997, zhan2019}.

\section{Sign Representation of Single-peaked Preferences} \label{sec:connect}
Let $\sigma$ be a linear order on $\mathcal{L}([n])$. 
A subset $X \subseteq[n] $ is called an \emph{ideal} of $\sigma$ if $x \in X$ and $y \prec_\sigma x$ imply $y \in X$.
In other words, if $\sigma=i_1 \cdots i_n$, then the ideal of $\sigma$ corresponds to initial segments $i_1 \cdots i_k$ $(k \le n)$ of linear order $\sigma$ and recursively, i.e., $i_1 \cdots i_j$ $(j \le k)$ \cite{DKK2012}.
A linear order or a permutation, $\sigma$, is \emph{consecutive} if all its ideals are ({consecutive}) intervals on $[n]$. 
There are totally $ \frac{n(n+1)}{2}+1$ intervals from length $0$ of trivial empty set to length $n$ on $[n]$.

Fix an alternative set $[n]$. 
By the definition given in Subsection \ref{sec:def_single}, each ideal of a single-peaked preference, or the set of the top to $p$th ($1 \le p \le n$) preferred  alternatives form a consecutive interval on $[n]$.  
Fig. \ref{fig:connected_4} is an example of all intervals for $n=4$, where each node represents an interval and arcs are drawn from the associated maximal ideals. 
From the definitions, the graphic structure presented in Fig \ref{fig:connected_4} is the same as the one in Fig. \ref{fig:single_tiling_4}, especially, there are at most two arcs from each node.
\setlength{\unitlength}{1mm}
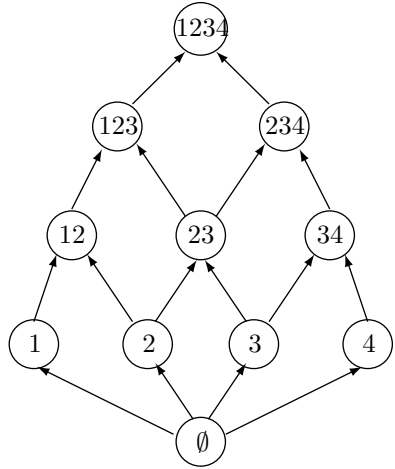
\begin{figure} 
 \centering
\begin{picture}(140,63)
\linethickness{0.6pt}

\put(37,-1){ \begin{minipage}[t]{62mm}
\begin{picture}(50,45)(0,0)
\linethickness{0.5pt}

\put(30,59){\circle{7}} \put(26.7,58){$1234$}

\put(21,49){\vector(1,1){7}} \put(39,49){\vector(-1,1){7}}

\put(41,46){\circle{6.5}} \put(38.5,45){$234$}
\put(19,46){\circle{6.5}} \put(16.5,45){$123$}
\put(32,34){\vector(2,3){6.5}}
\put(28,34){\vector(-2,3){6.5}}
\put(47,35){\vector(-1,2){4}}
\put(13,35){\vector(1,2){4}}

\put(13,31.5){\circle{6.5}} \put(11.3,30.5){$12$}
\put(30,31.5){\circle{6.5}} \put(28.3,30.5){$23$}
\put(47,31.5){\circle{6.5}} \put(45.3,30.5){$34$}

\put(52,20){\vector(-1,3){3}}
\put(39,20){\vector(2,3){6}}
\put(36,20){\vector(-2,3){5.5}}
\put(24,20){\vector(2,3){5.5}}
\put(21,20){\vector(-2,3){6}}
\put(8,20){\vector(1,3){3}}

\put(37,17){\circle{6.5}} \put(36.3,16){$3$}
\put(23,17){\circle{6.5}} \put(22.3,16){$2$}
\put(8,17){\circle{6.5}} \put(7.3,16){$1$}
\put(52,17){\circle{6.5}} \put(51.3,16){$4$}

\put(31,7){\vector(2,3){5}}
\put(29,7){\vector(-2,3){5}}
\put(26.4,5){\vector(-2,1){18}}
\put(33.3,5){\vector(2,1){18}}
\put(30,4){\circle{6.5}} \put(29.4,2.8){$\emptyset$}

\end{picture}
\end{minipage} }

\end{picture}
\caption{All intervals for $n=4$, directed paths represent the ideals of linear orders} 
\label{fig:connected_4}
\end{figure}

We now introduce a sign representation of single-peaked preferences.
\begin{theorem} \label{th:sign_exist}
Fix an alternative set $[n]$. Then, there exits a sequence of signs ``$+$'' and ``$-"$ expressions of length $n-1$ which uniquely specify each single-peaked preference.
\smallskip \\
\noindent 
{\rm {\sl Proof} \ 
We start from a trivial interval $\{i\} \subseteq [n]$ as the top alternative of preferences. 
We can obtain a single-peaked preference by adding an alternative in each step as follows until the interval is $[n]$. 
\begin{itemize}
\vspace{-1mm}
\item[] Let the current interval be $ \{j,j+1,\cdots,k \}$. 
For the next preferred alternative, there are at most two choices, $j-1$ (if $j-1 \ge 1$) or $k+1$ (if $k+1 \le n$), to keep the property of (consecutive) intervals.     
\end{itemize}
\vspace{-1mm}
The procedure means that each next preferred alternative in a single-peaked preference is uniquely determined by which neighbour is chosen based on the interval property. 
We denote the alternative by ``$+$'' if a larger number is added. Similarly, an alternative of the smaller number is represented by ``$-$''. 
The process terminates at $[n]$.
The number of ``$+$" or ``$-$" is uniquely determined by the top preferred alternative selected. 

It is obvious that all possible single-peaked preferences can be obtained from the construction, then we have the the claim.  
\pend }
\end{theorem} 

For example, single-peaked preferences $34251$, $43251$ with $n=5$ can also be expressed as $+-+-$ and $--+-$, respectively.
Conversely, $++-+$ means the single-peaked preference $23415$.\footnote{The author got this hint from more general cases of North-East (NE) lattice paths, see, e.g., \cite{Kratt2015}, where Catalan numbers \cite{abello_j1991} fit these structures.} 

 From the interval property of single-peaked preferences, the following combinatorial results can also be obtained easily \cite{BP2006}: 
\begin{itemize}
\item[-] the total number of preferences in $\mathcal{SP}([n])$ with $i \in [n]$ as the top preferred alternative is $\binom{n-1}{i-1}$,  here we assume $\binom{n-1}{0}=1$, 
\item[-]  $ \vert \mathcal{SP}([n]) \vert =
\sum_{i \in [n]} \binom{n-1}{i-1}=\binom{n-1}{0} + \binom{n-1}{1}+ \cdots + \binom{n-1}{n-1}=2^{n-1} $, 
\item[-] the numer $2^{n-1} $ can also be counted by choosing preferred alternatives individually such that the alternatives that are not chosen form the intervals 
(there are two choices, except for the last one).
\end{itemize}
 
Denote the position $p$ by $(p)$ in a single-peaked preference, where sign ``$+$'' or an increasing neighbour alternative is chosen. 
Then, single-peaked preferences $34251$, $43251$ with $n=5$ can be written as $(2)(4)$ and $(4)$, respectively. 
Conversely, $(2)(3)(5)$ means a single-peaked preference $23415$.    
The position $\bar{p}$ denoted by $({\bar p})$ where the sign ``$-$'' or a decreasing alternative is chosen can be treated symmetrically.
Call such $(p)$ and $(\bar{p})$ the \emph{positive position} and  \emph{negative position}, respectively. 
Let $\sigma$ be a single-peaked order, and denote the set of  all negative positions of $\sigma$ by Neg($\sigma$).
Now we have following Proposition \ref{th:bruhat_num}.
\begin{proposition} \label{th:bruhat_num}
For each $\sigma \in \mathcal{SP}([n]),$ we have $ \vert{\rm Inv}(\sigma)\vert = \sum_{(\bar{p})\in {\rm Neg}(\sigma)} (\bar{p}-1) $. 
\smallskip \\
\noindent 
{\rm {\sl Proof} \ Based on the interval property of single-peaked preferences, each $(\bar{p})$ means an alternative $i$ which is the minimum natural number from the top to the $\bar{p}$th preference in $\sigma$. 
From the proof of Theorem \ref{th:sign_exist}, appending $i$ of $(\bar{p})$ adds exactly $\bar{p}-1$ inversion pairs to $ {\rm Inv}(\sigma)$. 
While positive position $(p)$ is the maximal natural number from the top to the $p$th preference in $\sigma$, these alternatives  contribute nothing to $ {\rm Inv}(\sigma)$. 
Since the above discussion applies to any alternatives of two cases ``$-$" and ``$+$" in $\sigma$, we have the conclusion. 
\pend }
\end{proposition} 
For a special case of Proposition \ref{th:bruhat_num}, the reverse linear order $n \cdots 1$, we have $1+\cdots + (n-1)= n(n-1)/2=\binom{n}{2}$. 

The following basic operations defined on the sign representation  characterise the associated Bruhat poset of each maximal single-peaked domain. 
\begin{theorem} \label{th:con_sign}
Fix a single-peaked domain $ \mathcal{SP}([n]) $. 
The following (i) and (ii) are the only operations traversing along the edges of associated Bruhat poset $\mathbb{B}(\mathcal{SP}[n],1)$ based on the sign representation:
\vspace{1mm}

(i) change the first sign, 

(ii) swap two opposite neighbour signs, from ``$+-$" to ``$-+$", or from ``$-+$" to ``$+-$".  
\smallskip \\
\noindent 
{\rm {\sl Proof} \ 
The above operation (i) means that we swap the top and second preferred alternatives; the resulting linear order is single-peaked since the two alternatives form an interval. 
The operation (ii) means: we select the opposite neighbour, the small number, instead of the current larger one, and vice versa, in creating single-peaked preferences in the proof of Theorem  \ref{th:sign_exist}. 
We see, also by Proposition \ref{th:bruhat_num}, that both operations (i) and (ii), create or reduce, the number of inversion pairs by one.  
Flipping the top sign changes the top-preferred alternative, hence the total number of positive or negative signs.  
  
For the ``only" part of the claim, if we change the middle sign or swap the signs that are not neighbours, by Proposition \ref{th:bruhat_num}, the difference of the number or the cardinality of inversion pairs will be larger than 1.
It should be pointed out that the sign position is smaller than the associated alternative position by one since the sign represents the comparison with the previous neighbour alternatives.   
\pend }
\end{theorem}

From the proof of Theorem \ref{th:sign_exist}, we know that two alternatives or two natural numbers associated with two consecutive same sign, ``$++$" or ``$--$", are forbidden crossing to keep the property of intervals, otherwise, they would bring no chang or cause a failure in the sign representation.
This is an essential fact. 
Therefore, from Theorem \ref{th:sign_exist} and Theorem \ref{th:con_sign},
the introduced sign representation specifies \emph{exactly} both the basic structure and operations of single-peaked preferences.

In Theorem \ref{th:con_sign}, when two more operations are simultaneously possible in a current sign representation, we obtain single-peaked linear orders that belong to an \emph{equivalent class} \cite{g_reiner08}.

As an application of the interval property and sign representation, we show the characterization of single-peaked preferences by Puppe \cite{EFS2020, puppe2018}. 

A domain $\mathcal{D} \subset \mathcal{L}([n])$ is called \emph{minimally rich} if each alternative is top preferred in at least one linear order of $\mathcal{D}$.
A domain $\mathcal{D}$ is \emph{semi-connected} if it contains two completely reversed orders, and a directed path connecting them in the digraph of Bruhat poset $\mathbb{B}(\mathcal{D},1)$ (here, the connectedness of the digraph of Bruhat poset is defined naturally). 
Without loss of generality, we assume that a single-peaked domain contains the following two reverse orders, $\alpha \equiv 12\cdots (n-1)n$ and $\omega \equiv n(n-1)\cdots 21$, with Inv$(\alpha)=\emptyset$ and Inv$(\omega)=\Omega$. 
A domain that contains such a pair of reverse orders is said to
have \emph{maximal width}. 
Each single-peaked domain (including more general case) can contain at most one pair of such reverse orders \cite{slinko2019}. 
The semi-connected condition means that there are at least $\binom{n}{2}+1$ linear orders on a domain $\mathcal{D} \subset \mathcal{L}([n])$ \cite{EFS2020, puppe2018}. 

Here is a different approach to show partly a known result.   

\begin{theorem}[Theorem 1, \cite{puppe2018}] \label{th:min_rich}
The domain $\mathcal{SP}(n])$ of all single-peaked orders is a connected and minimally rich CD of maximal width. 
In particular, $\mathcal{SP}([n])$ is semi-connected.

\smallskip 
{\rm
Since single-peaked linear orders can be constructed starting from any alternative, as shown in the proof of Theorem \ref{th:sign_exist}, which proves the minimal richness. 
Two reverse orders are constructed by starting from ``$1$" to obtain $\alpha$, or ``$n$" to get a reverse liner order $\omega$, hence $\mathcal{SP}(n])$ is also maximal width. 

For the semi-connected condition, we begin from $\alpha \equiv 12\cdots (n-1)n$. 
If there is a sign ``$-$" before a neighbour sign ``$+$", by Theorem \ref{th:sign_exist} and Theorem \ref{th:con_sign}, operation (ii) along a directed arc of associated Bruhat poset $\mathbb{B}(\mathcal{D},1)$ leads to a single-peaked preference with one more reverse pair.  
Otherwise, we have $+\cdots+-\cdots -$; the latter part may be empty.
Then, by operation (i) of Theorem \ref{th:con_sign}, we can change the first sign ``$+$" to ``$-$" and keep the single-peaked linear orders. 
Starting from $\alpha$, or $+ \cdots +$, by the above arguments, 
we can reach the reverse $\omega \equiv n \cdots 1$ with $- \cdots -$.
Both operations (i) and (ii) are associated with an arc of $\mathbb{B}(\mathcal{D},1)$ according to the results of Theorem \ref{th:con_sign}.
This proves the semi-connectedness.
\pend 
}
\end{theorem}
\section{Concluding Remarks} \label{sec:conclud}

The contribution of this study related to the interval property and sign representation of single-peaked preferences are summarized as follows.
\begin{description}
\item[Sign representation of single-peaked preferences]\mbox{}\\
\begin{itemize}
\vspace{-4mm}
\item[-] Although single-peaked preferences are well studied since 1948 \cite{black1948}, to the best of our knowledge, the sign representation of single-peaked preferences introduced here that procedurally and explicitly contains inversion information is new. 
\item[-] Two basic operations are defined based on the sign representation  which completely define the structure of Bruhat poset induced by inversion pair inclusion. 
We also provide some applications. 
\end{itemize}
\end{description}

We hope that this note provides some intuitions on the single-peaked preferences and related topics for further studies. 

\Delete{
The ideals of preferences share the similar property (axiom) of {feasible ordering} of the greedoids \cite{s2021}. 
Greedy algorithms work on a rooted-tree structure. 
Trees are more general single-peaked domains \cite{d1982}.

A preference domain is \emph{single-crossing} if voters (or preferences) can be ordered as $(v_1, \dots, v_m)$ such that for each pair $\{i,j\}$ of alternatives, there exists $k \in [m]$ with $\{ \ell \in  [m] \mid i \succ_\ell j \} = [k] $ \cite{sww2021}. 
Restricted on single-crossingness, a \emph{representative} voter of majority \cite{rothst1991}, or a solution on the greedoid of associated weight \cite{s2021}, exists.
}

\medskip

\bmhead{Acknowledgments}
The author's work was partially supported by JSPS KAKENHI Grant Number 20K04970.

\medskip
\noindent
{\small {\bf Data Availability}  The datasets generated during the current study are available from the author on
reasonable request.}

\subsection*{Conflict of interest}
The author declares that there is no conflict of interest.




\end{document}